\documentclass{article}
\usepackage{cite}
\usepackage{amsmath,amssymb,amsfonts}
\usepackage{algorithmic}
\usepackage{graphicx}
\usepackage{textcomp}

\usepackage{epsfig}
\usepackage{latexsym}
\usepackage{pdfsync}

\title{ Boundary Control of the Wave Equation \\via Linear Quadatic Regulation\thanks{Submitted to the editors DATE.
\funding{This work was funded by AFOSR under under grant XXXXX }}}

\author{Arthur J. Krener \thanks{Department of Applied Mathematics, Naval Postgraduate School, Monterey CA
  (\email{ajkrener@nps.edu})}}

\newcommand{\bc}{\begin{center}}
\newcommand{\ec}{\end{center}}

\newcommand{\eq}{\begin{equation}\begin{array}{rllllllllllllllllllllllllllllllll}}
\newcommand{\ee}{\end{array}\end{equation}}
\newcommand{\bmt}{\left[ \begin{array}{ccccccccc}}
\newcommand{\emt}{\end{array}\right]}
\newcommand{\bea}{\begin{eqnarray}}
\newcommand{\eea}{\end{eqnarray}}
\newcommand{\bean}{\begin{eqnarray*}}
\newcommand{\eean}{\end{eqnarray*}}

\begin{document}

\Large 
\bc

Boundary Control of the Wave Equation \\via Linear Quadatic Regulation\\
\vspace{0.2in}
Arthur J. Krener\\ Department of Applied Mathematics\\ Naval Postgraduate School\\
Monterey, CA, USA
\\
ajkrener@nps.edu
\ec
\normalsize

\begin{abstract}
  We consider the Linear Quadratic Regulation  for the  boundary  control 
  of the one dimensional linear wave equation under both Dirichlet
  and Neumann activation.  For each activation we present
  a Riccati partial differential equation that we explicitly solve.
  The derivation the Riccati partial differential equations is by the
  simple and explicit technique of completing the square.
  \end{abstract}

\section{Introduction}
The control of infinite dimensional systems is well treated in numerous works, \cite{JL71}, \cite{CZ95}, \cite{LT00}.
   \cite{KS08}, \cite{CZ20}.   In particular Lasiecka and Triggiani \cite{LT86} prove the existence of solutions to the  Riccati
   equations that arise in the boundary control of  hyperbolic problems in  higher dimensions.  
   Burns and King \cite{BK95} obtain integral representations for the feedback operators for hyperbolic problems with 
   Kelvin-Voight damping and non-compact input operators.

   More recently \cite{Kr20} we were able to explicitly solve a boundary control problem for a nonlinear reaction diffusion 
   equation by completing the square and extending Al'brekht's method \cite{Al61}. 
In this paper we explicitly find infinite horizon Linear Quadratic Regulator (LQR) for the one dimensional 
linear wave equation  under either Dirichlet and Neumann activation.   We use the simple and
explicit technique of completing the square to derive Riccati partial differential equations
for the optimal cost and optimal control.  By decoupling the spatial frequencies we reduce
Riccati partial differential equations to an infinite family of two dimensional algebraic Riccati
equations that can be solved explicitly.

 \section{Dirichlet Boundary Control of the Linear Wave Equation via LQR}
Consider the linear wave equation subject to Dirichlet  boundary control at one end.
We could consider  Dirichlet  boundary control at both ends but to keep the 
notation simple we do not do so.
The spatial variable is $x\in[0,1]$ and the model is
\bean
\frac{\partial^2 w}{\partial t^2}(x,t)&=&\frac{\partial^2 w}{\partial x^2}(x,t)-\alpha\frac{\partial w}{\partial t}(x,t) \\
w(0,t)=\beta  u(t),&&w(1,t)=0\\
w(x,0)=w_0(x),&& \frac{\partial w}{\partial t}(x,0)=w_1(x,0)
\eean
where $\alpha$ is a nonegative constant.   If $\alpha>0$ the wave equation is damped
and if $\alpha=0$ it is undamped.  For Dirichlet activation is convenient to have the control act at $x=0$.

Suppose $w^*(x,t)$ is  any solution of unforced wave equation subject to zero Dirichlet boundary conditions and $w(x,t)$ is a solution 
of wave equation subject to Dirichlet boundary control.  Then we define
\bea \label{z}
z(x,t)&=& \bmt z_1(x,t)\\z_2(x,t)\emt \ = \ \bmt w(x,t)-w^*(x,t)\\  \frac{\partial w}{\partial t}(x,t)- \frac{\partial w^*}{\partial t}(x,t) \emt
 \eea
and $z(x,t)$ satisfies
\bean \label{vpde}
\frac{\partial z}{\partial t}(x,t)&=& {\cal A} z(x,t)
 \eean 
 subject to boundary and initial conditions
 \bean
  z_1(0,t)=\beta u(t),&\quad\ & z_1(0,t)=0
  \\
 z_1(x,0)=w_0(x)-w^*(x,0),&\quad\ &z_2(x,0)=w_1(x,0) - \frac{\partial w^*}{\partial t}(x,0)
 \eean
  where $ {\cal A}$ is the matrix differential operator
 \bea \label{mdo}
 {\cal A}&=& \bmt 0& 1\\\frac{\partial^2 }{\partial x^2}&-\alpha \emt
 \eea
 
 The open loop eigenvalues of ${\cal A}$  under Dirichlet boundary conditions are 
 \bea \label{laD}
  \lambda_n &=&-{\alpha+\mbox{ sign}(n) \sqrt{\alpha^2-4 n^2\pi^2}\over 2}
  \eea
  for $n=\pm1,\pm 2,\ldots $.
The coresponding  eigenvectors are
\bean
\phi^n(x)=\bmt {1\over \lambda_n} \\1\emt \sin |n|\pi x
\eean 
 
Notice that if the equation is undamped, $\alpha=0$, then all the eigenvalues
are imaginary.  If $\alpha>0$ then some of the eigenvalues might be negative
real numbers but as $|n|$ gets larger they become complex numbers with negative real part.

We wish to stabilize this system to $z(x,t)=0$
so we set up a LQR  problem.
We choose a $2\times  2$ dimensional matrix valued function
\bean
Q(x_1,x_2)&=& \bmt Q_{1,1}(x_1,x_2) &Q_{1,2}(x_1,x_2)\\Q_{2,1}(x_1,x_2)&Q_{2,2}(x_1,x_2)\emt
\eean
which is symmetric and nonnegative definite for each $(x_1,x_2)\in {\cal S}= [0,1]\times [0,1]$
and symmetric with respect to $(x_1,x_2)$,
$
Q(x_1,x_2)=Q(x_2,x_1)
$.
We also choose a positive scalar $R$
and consider the problem of minimizing by choice
of control $u(t)$ the criterion
\bea \label{crit}
\int_0^\infty \iint_{\cal S}  z'(x_1,t) Q(x_1,x_2) z(x_2,t)\ dA+R( u(t))^2\ dt
\eea
where $dA=dx_1dx_2$.

Let $P(x_1,x_2)$ be a  $2\times  2$ dimensional symmetric matrix valued function
\bean
P(x_1,x_2)&=& \bmt P_{1,1}(x_1,x_2) &P_{1,2}(x_1,x_2)\\P_{2,1}(x_1,x_2)&P_{2,2}(x_1,x_2)\emt
\eean
 which is also symmetric in $(x_1,x_2)$,
$
P(x_1,x_2)=P(x_2,x_1)
$.
Suppose there exists a control trajectory $u(t)$ such that the resulting state trajectory $z(x,t) \to 0$ as $t \to \infty$.  Then by the Fundamental 
Theorem of Calculus
\bean
0&=&\iint_{\cal S}  z'(x_1,0) P(x_1,x_2) z(x_2,0)\ dA\\
&& +\int_0^\infty \iint_{\cal S} {d\over dt}\left( z'(x_1,t) P(x_1,x_2) z(x_2,t)\right)\ dA \ dt
\eean

We expand the time integrand   into components

\bean
0&=&\iint_{\cal S}  z'(x_1,0) P(x_1,x_2) z(x_2,0)\ dA
\\
&& +\int_0^\infty \iint_{\cal S} z_2(x_1,t)P_{1,1}(x_1,x_2)z_1(x_2,t)+ z_1(x_1,t)P_{1,1}(x_1,x_2)z_2(x_2,t)\\
&&+z_2(x_1,t)P_{1,2}(x_1,x_2)z_2(x_2,t)+ z_2(x_1,t)P_{2,1}(x_1,x_2)z_2(x_2,t)\\
&&+\left(  \frac{\partial^2 z_1}{\partial x_1^2}(x_1,t)-\alpha z_2(x_1,t)\right)P_{2,1}(x_1,x_2)z_1(x_2,t)\\
&&+z_1(x_1,t)P_{1,2}(x_1,x_2)\left(\frac{\partial^2 z_1}{\partial x_2^2}(x_2,t)-\alpha z_2(x_2,t)\right)\\
&&+ \left( \frac{\partial^2 z_1}{\partial x_1^2}(x,t)-\alpha z_2(x_1,t)\right)P_{2,2}(x_1,x_2)z_2(x_2,t)\\
&&+z_2(x_1,t)P_{2,2}(x_1,x_2)\left(\frac{\partial^2 z_1}{\partial x_2^2}(x,t)-\alpha z_2(x_2,t)\right)
\ dA \ dt
\eean

We assume that $P(x_1,x_2)$ also satisfies Dirichlet boundary conditions
\bean
P(0,x_2)&=& P(0,x_2)\ =\ 0\\
P(x_1,0)&=&P(x_1,0)\ =\ 0
\eean
 and we integrate by parts twice with respect to the $x_i$ to get

\bean
0&=&\iint_{\cal S}  z'(x_1,0) P(x_1,x_2) z(x_2,0)\ dA\\
&& +\int_0^\infty \iint_{\cal S} z_2(x_1,t)P_{1,1}(x_1,x_2)z_1(x_2,t)+ z_1(x_1,t)P_{1,1}(x_1,x_2)z_2(x_2,t)\\
&&+z_2(x_1,t)P_{1,2}(x_1,x_2)z_2(x_2,t)+ z_2(x_1,t)P_{2,1}(x_1,x_2)z_2(x_2)\\
\\&& -\alpha z_2(x_1,t)P_{2,1}(x_1,x_2)z_1(x_2)-\alpha z_1(x_1,t)P_{1,2}(x_1,x_2)z_2(x_2)
\\ && -2\alpha z_2(x_1,t)P_{2,2}(x_1,x_2)z_2(x_2)\\
&&+  z_1(x_1,t)\frac{\partial^2 P_{2,1}}{\partial x_1^2}(x_1,x_2)z_1(x_2,t)+z_1(x_1,t)\frac{\partial^2 P_{1,2}}{\partial x_2^2}(x_1,x_2)z_1(x_2,t)\\
&&+  z_1(x_1,t)\frac{\partial^2 P_{2,2}}{\partial x_1^2}(x_1,x_2)z_2(x_2,t)+z_2(x_1,t)\frac{\partial^2 P_{2,2}}{\partial x_2^2}(x_1,x_2)z_1(x_2,t)
\ dA \ dt \\
&&+\int_0^\infty  \int_0^1 \beta u(t)\frac{\partial P_{2,1}}{\partial x_1} (0,x_2) z_1(x_2,t)\ dx_2\ dt
\\
&&+\int_0^\infty  \int_0^1 z_1(x_1,t)\frac{\partial P_{1,2}}{\partial x_2} (x_1,0)\beta u(t) \ dx_1 \ dt
\\
&&+\int_0^\infty  \int_0^1 \beta u(t)\frac{\partial P_{2,2}}{\partial x_1} (0,x_2) z_2(x_2,t)\ dx_2\ dt
\\
&&+\int_0^\infty  \int_0^1 z_2(x_1,t)\frac{\partial P_{2,2}}{\partial x_2} (x_1,0)\beta u(t) \ dx_1
 \ dt
\eean

We add this to the criterion (\ref{crit}) to get the equivalent criterion to be minimized

\bean
&&\iint_{\cal S}  z'(x_1,0) P(x_1,x_2) z(x_2,0)\ dA
\\
&& +\int_0^\infty \iint_{\cal S} z_2(x_1,t)P_{1,1}(x_1,x_2)z_1(x_2,t)+ z_1(x_1,t)P_{1,1}(x_1,x_2)z_2(x_2,t)\\
&&+z_2(x_1,t)P_{1,2}(x_1,x_2)z_2(x_2,t)+ z_2(x_1,t)P_{2,1}(x_1,x_2)z_2(x_2)\\
&& -\alpha z_2(x_1,t)P_{2,1}(x_1,x_2)z_1(x_2)-\alpha z_1(x_1,t)P_{1,2}(x_1,x_2)z_2(x_2)
\\ && -2\alpha z_2(x_1,t)P_{2,2}(x_1,x_2)z_2(x_2)\\
&&+  z_1(x_1,t)\frac{\partial^2 P_{2,1}}{\partial x_1^2}(x_1,x_2)z_1(x_2,t)+z_1(x_1,t)\frac{\partial^2 P_{1,2}}{\partial x_2^2}(x_1,x_2)z_1(x_2,t)\\
&&+  z_1(x_1,t)\frac{\partial^2 P_{2,2}}{\partial x_1^2}(x_1,x_2)z_2(x_2,t)+z_2(x_1,t)\frac{\partial^2 P_{2,2}}{\partial x_2^2}(x_1,x_2)z_1(x_2,t)\\
&&+z'(x_1,t) Q(x_1,x_2)z(x_2,t)\ dA +R(u(t)^2\ dt 
\\&&+\int_0^\infty  \int_0^1 \beta u(t)\frac{\partial P_{2,1}}{\partial x_1} (0,x_2) z_1(x_2,t)\ dx_2\ dt
\\
&&+\int_0^\infty  \int_0^1 z_1(x_1,t)\frac{\partial P_{1,2}}{\partial x_2} (x_1,0)\beta u(t) \ dx_1\ dt
\\
&&+\int_0^\infty  \int_0^1 \beta u(t)\frac{\partial P_{2,2}}{\partial x_1} (0,x_2) z_2(x_2,t)\ dx_2\ dt
\\
&&+\int_0^\infty  \int_0^1 z_2(x_1,t)\frac{\partial P_{2,2}}{\partial x_2} (x_1,0)\beta u(t) \ dx_1
 \ dt
\eean

We would like to choose 
\bean
K(x) &=& \bmt K_{1}(x) & K_{2}(x) \emt
\eean
so that the sum of the integrands of the time integrals is  a perfect square of the form
\bean
\iint_{\cal S} \left(u(t)-K(x_1) z(x_1,t)\right)'R\left(u(t)-K(x_2) z(x_2,t)\right)\ dA
\eean 
Obviously the  terms quadratic in $u(t)$  match  so we equate terms containing the product of  $u(t)$ and $ z(x_2,t)$.   This leads to the equation
\bea \label{KKK}
-R K(x_2) &=&\beta \bmt
\frac{\partial P_{2,1}}{\partial x_1} (0,x_2) &\frac{\partial P_{2,2}}{\partial x_1} (0,x_2) \emt 
\eea
so we assume
\bea \label{K1}
K(x_2)&=&\bmt K_1(x_2)&K_2(x_2)\emt \ =\  -R^{-1} \beta  \bmt \frac{\partial P_{2,1}}{\partial x_1} (0,x_2) & \frac{\partial P_{2,2}}{\partial x_1} (0,x_2) \emt 
\eea
By symmetry
\bean
K(x_1)&=&\bmt K_1(x_1)&K_2(x_1)\emt \ =\  -R^{-1} \beta  \bmt \frac{\partial P_{1,2}}{\partial x_1} (x_1,0) & \frac{\partial P_{2,2}}{\partial x_1} (x_1,0) \emt 
\eean

Next we equate terms quadratic in $z(x,t)$ and we obtain the equations
\bea
\nonumber&&\frac{\partial^2 P_{2,1}}{\partial x_1^2}(x_1,x_2) + \frac{\partial^2 P_{1,2}}{\partial x_2^2}(x_1,x_2)+Q_{1,1}(x_1,x_2)
\\&&=\gamma^2   \frac{\partial P_{1,2}}{\partial x_2} (x_1,0) \frac{\partial P_{2,1}}{\partial x_1} (0,x_2)
\label{RPDE11a}\\
\nonumber&&P_{11}(x_1,x_2)-\alpha P_{1,2}(x_1,x_2) +\frac{\partial^2 P_{2,2}}{\partial x_1^2}(x_1,x_2)+Q_{1,2}(x_1,x_2)   \\
&&   =\gamma^2   \frac{\partial P_{1,2}}{\partial x_2} (x_1,0) \frac{\partial P_{2,2}}{\partial x_1} (0,x_2)
\label{RPDE12a}\\
\nonumber&&P_{11}(x_1,x_2)-\alpha P_{2,1}(x_1,x_2)+\frac{\partial^2 P_{2,2}}{\partial x_2^2}(x_1,x_2)+Q_{21}(x_1,x_2) 
\\ && =\gamma^2\frac{\partial P_{2,2}}{\partial x_2} (x_1,0)\frac{\partial P_{2,1}}{\partial x_1} (0,x_2) \label{RPDE21a}\\
\nonumber&&P_{1,2}(x_1,x_2)+P_{2,1}(x_1,x_2)-2\alpha P_{2,2}(x_1,x_2)+Q_{22}(x_1,x_2)\\
&&  =\gamma^2 \frac{\partial P_{2,2}}{\partial x_2} (x_1,0)\frac{\partial P_{2,2}}{\partial x_1} (0,x_2) \label{RPDE22a}
\eea
where $\gamma^2=R^{-1} \beta^2$.

We call these equations (\ref{RPDE11a},\ref{RPDE12a},\ref{RPDE21a},\ref{RPDE22a}) the Riccati PDE 
for the Dirichlet LQR control of the wave equation.

We assume that $P(x_1,x_2)$ can be expressed in terms of the eigenfunctions of $\frac{\partial^2}{\partial x^2}$
with respect to Dirichlet  boundary conditions  at $x=0$ and $x=1$.
\bea  \label{Psum}
P(x_1,x_2)&=& \sum_{m=1}^\infty \sum_{n=1}^\infty \bmt P^{m,n}_{1,1}&P^{m,n}_{1,2}\\P^{m,n}_{2,1}&P^{m,n}_{2,2}\emt \sin m\pi x_1  \sin n\pi x_2
\eea
and  $P^{m,n}_{1,2}=P^{m,n}_{2,1}$.
Then (\ref{K1}) implies 
\bea \label{Ksum}
K(x_2)&=&-R^{-1}\beta  \sum_{m=1}^\infty  \sum_{n=1}^\infty \bmt P_{2,1}^{m,n} &P_{2,2}^{m,n} \emt m\pi \sin n\pi x 
\eea
 
We assume that $Q(x_1,x_2)$ has a similar expansion
\bea \label{Qsum}
Q(x_1,x_2)&=& \sum_{m=1}^\infty \sum_{n=1}^\infty  Q^{m,n} \sin m\pi x_1  \sin n\pi x_2
\eea
where
\bean
Q^{m,n}&=& \bmt Q^{m,n}_{1,1}&Q^{m,n}_{1,2}\\Q^{m,n}_{2,1}&Q^{m,n}_{2,2}\emt 
\eean
and  $Q^{m,n}_{1,2}=Q^{m,n}_{2,1}$.

This leads to an infinite algebraic Riccati equation for $P^{m,n}$.
To
 simplify the analysis we decouple the spatial frequencies   by  assuming
 \bea \label{Q}
Q^{m,n}&=& \delta_{m,n}\bmt Q^{n,n}_{1,1}&Q^{n,n}_{1,2}\\Q^{n,n}_{2,1}&Q^{n,n}_{2,2}\emt 
\eea

Then we guess 
\bean
 P^{m,n}&=& \delta_{m,n}  P^{n,n}
\eean

If this is true then the Riccati PDE   (\ref{RPDE11a},\ref{RPDE12a},\ref{RPDE21a},\ref{RPDE22a}) implies that 
\bea
0&=&   -2n^2\pi^2 P^{n,n}_{1,2} +Q^{n,n}_{1,1}-n^2\pi^2\gamma^2\left(P^{n,n}_{1,2}\right)^2  \label{11}\\
0&=& P^{n,n}_{1,1}-\alpha P^{n,n}_{1,2} -n^2\pi^2 P^{n,n}_{2,2}+Q^{n,n}_{1,2} -n^2\pi^2\gamma^2P^{n,n}_{1,2}P^{n,n}_{2,2} \label{12}\\
0&=& P^{n,n}_{1,1} -\alpha P^{n,n}_{2,1}-n^2\pi^2 P^{n,n}_{2,2}+Q^{n,n}_{2,1} -n^2\pi^2\gamma^2P^{n,n}_{2,2}P^{n,n}_{2,1}\label{21}\\
0&=& 2P^{n,n}_{1,2}-2\alpha P^{n,n}_{2,2}+Q^{n,n}_{2,2}-n^2\pi^2 \gamma^2\left(P^{n,n}_{2,2}\right)^2 \label{22}
\eea
where $\gamma^2= R^{-1}\beta^2$.

For each $n=1,2,\ldots$ these are the Riccati equations of  the two dimensional  LQR with matrices 
\eq  \label{lqr2}
F^{n,n}=\bmt 0&1\\ -n^2 \pi^2 &-\alpha\emt, && G^{n,n}=\bmt 0\\n\pi \beta\emt\\
Q^{n,n}=\bmt Q^{n,n}_{1,1}& Q^{n,n}_{1,2}\\Q^{n,n}_{2,1}&Q^{n,n}_{2,2}\emt,&& R^{n,n}=\bmt R\emt
\ee

We use the quadratic formula to solve   (\ref{11})
\bea \label{p12}
P^{n,n}_{1,2}&=& {-1\pm \sqrt{1+{\gamma^2Q^{n,n}_{1,1}\over n^2\pi^2}}\over \gamma^2}
\eea
then the quadratic formula applied to (\ref{22})  implies
\bea \label{p22}
P^{n,n}_{2,2}&=&{-\alpha \pm \sqrt{\alpha^2+n^2\pi^2\gamma^2\left(Q^{n,n}_{2,2}+2P^{n,n}_{1,2}\right)}\over n^2\pi^2\gamma^2}
\eea
Since we want $P^{n,n}_{2,2}$ to be nonnegative we need to take the positive sign.
Then  (\ref{12}) implies 
\bea \label{p11}
P^{n,n}_{1,1}&=& \alpha P^{n,n}_{1,2}+n^2\pi^2 \left(1+\gamma^2 P^{n,n}_{1,2}  \right) P^{n,n}_{2,2}-Q^{n,n}_{1,2}
\eea
%
If the two dimensional LQR (\ref{lqr2}) satisfies the standard conditions then the associated Riccati equation has a unique nonnegative definite solution.  This implies that if we take the negative sign in (\ref{p12}) the resulting $P^{n,n}$ is not nonnegative definite.

The $2\times 2$ closed loop system is 
\bean
F^{n,n}+GK^{n,n}&=& \bmt 0&1\\ -n^2\pi^2 -\gamma^2 P_{2,1}^{n,n}&-\alpha -\gamma^2 P_{2,2}^{n,n}\emt
\eean
and the closed loop eigenvalues are
\bea \nonumber
\mu_n&=& 
-{\alpha +n^2\pi^2\gamma^2P^{n,n}_{2,2}\over 2}+\left(\mbox{sign n}\right){\sqrt{(\alpha +n^2\pi^2 \gamma^2P^{n,n}_{2,2})^2-4(n^2\pi^2+\gamma^2 P^{n,n}_{2,1})}\over 2}\\  \label{cleva}
\eea
for $n=\pm 1,\pm2, \pm3,\ldots$.
For $n=1,2,\ldots$ the corresponding eigenvectors of the $n^{th}$ $2\times 2$ closed loop system are 
 \bea \label{eigv2}
 \bmt {1\over \mu_n}\\1 \emt
 ,&\quad\ & \bmt {1\over \mu_{-n}}\\1 \emt
 \eea
 The corresponding eigenvectors of the infinite closed loop system are 
 \bea \label{eigv2in}
v_n(x)= \bmt {1\over \mu_n}\\ 1\emt \sin |n|\pi x  ,&\quad\ &  v_{-n}(x)= \bmt {1\over \mu_{-n}}\\ 1\emt \sin |n|\pi x
 \eea
 Notice $\mu_n$ and $\mu_{-n}$ are complex conjugates as are $v_n(x)$ and $v_{-n}(x)$.
 
The trajectories of the infinite closed loop system are 
\bean
z(x,t)&=& \sum_{n=-\infty}^\infty \zeta_n(t) \bmt {1\over \mu_n}\\1 \emt \sin |n|\pi x
\eean
where
\bean
\zeta_n(t) &=& e^{\mu_n t}  \zeta_n^0
\eean
Since we want a real valued $z(x,t)$,  $\zeta_n^0$ and $\zeta_{-n}^0$ must be complex conjugates.

The  quadratic kernel of the optimal cost is given by (\ref{Psum}) which reduces to 
\bea \label{Psum_n}
P(x_1,x_2)&=& \sum_{n=1}^\infty \bmt P_{1,1}^{n,n}&P_{1,2}^{n,n}\\P_{2,1}^{n,n}&P_{2,2}^{n,n}\emt \sin n\pi x_1 \sin n\pi x_2
\eea
Notice we can control each spatial frequency indepedently.   If we don't want to damp out the $n^{th} $ spatial frequency then we set
$Q^{n,n}=0$ so that  $P^{n,n}=0$ and $K^{n,n}=0$. 

A critical issue is whether the series (\ref{Psum_n}) is convergent.  Clearly some norm of $Q^{n,n}$ needs to go to zero faster than ${1\over n}$ for the series
(\ref{Qsum}) to converge.  But probably $Q^{n,n}$ needs to go to zero even faster for (\ref{Psum_n}) to converge.  Of course if the system is damped we can set
$Q^{n,n}=0$ for $n>N$ and let the damping stabilize the higher spatial modes.

But  if the system is undamped
  how many spatial modes can we stabilize and how fast can we dampen them?  To answer these questions we look at a simple example,  we assume 
\bean 
\beta=1\\
R=1\\
Q^{n,n}_{1,1}=Q^{n,n}_{2,2}={q\over n^r }\\
Q^{n,n}_{1,2}=Q^{n,n}_{2,1}=0
\eean
For (\ref{Qsum})  to converge we must take $r>1$. 


For this  system
\bean
\gamma^2&=&1\\
P^{n,n}_{1,2}&=&-1+ \sqrt{1+{q\over \pi^2 n^{2+r}}}\\
P^{n,n}_{2,2}&=&\sqrt{{q/n^r+2P^{n,n}_{1,2}\over n^2\pi^2}}
\\
P^{n,n}_{1,1}&=& n^2\pi^2\left(1+P^{n,n}_{1,2}\right) P^{n,n}_{2,2}
\eean

 First note that 
by the Mean Value Theorem there exists an $s$ between $ 1$ and  $1+{q\over \pi^2n^ {2+r}}$
such that
\bean
P^{n,n}_{1,2}&=& \sqrt{1+{q\over \pi^2 n^{r+2}}}-1\ =\ {1\over 2s^{1/2}} {q\over \pi^2 n^{2+r}}
\eean 
The function ${1\over 2s^{1/2}}$ is montonically decreasing between $ 1$ and  $1+{q\over \pi^2n^{2+r}}$
and takes on its maximum value ${1\over 2}$ at 	$s= 1$ so
\bean
0&<&P^{n,n}_{1,2}\le {1\over 2}  {q\over \pi^2 n^{r+2}}
\eean
Then
\bean
P^{n,n}_{2,2}&=&\sqrt{{q/n^r+2P^{n,n}_{1,2}\over n^2\pi^2}} \ \le\ { c_1 \over n^{1+r/2}}
\\
P^{n,n}_{1,1}& \le & { c_2 \over n^{-1+r/2}}
\eean
for some positive constants $c_1,c_2$.  For  $P_{1,1}(x_1,x_2)$ term in (\ref{Psum_n}) to converge  we must take $r>4$. 
Then $P^{n,n}_{2,2}< {c_1\over n^s}$ for some $s>3$.  If we choose $r$ a little bigger than $4$ then the term outside the square root in the closed loop eigenvalues (\ref{cleva}), $-{n^2\pi^2\gamma^2P^{n,n}_{2,2}\over 2}$ is decaying so the higher spatial modes are substatially less damped.   Moreover the first term inside the the square root in (\ref{cleva}) is going to zero while the negative term is   growing so at least the higher eigenvalues are complex.

But for the optimal feedback (\ref{Ksum}) to converge  we only need to take $r>2$. If we take $r$ a little larger than $2$
then the term outside the square root in the closed loop eigenvalues (\ref{cleva}), $-{n^2\pi^2\gamma^2P^{n,n}_{2,2}\over 2}$, is decaying like $n^{-s}$ for some $s$ a little larger than $0$.   So again  the higher  the spatial mode the less damping.

 \section{Neumann Boundary Control of the Linear Wave Equation via LQR}
   In this section we   assume Neumann boundary control instead of Dirichlet boundary control.
   The model is now
\bean
\frac{\partial^2 w}{\partial t^2}(x,t)&=&\frac{\partial^2 w}{\partial x^2}(x,t)-\alpha\frac{\partial w}{\partial t}(x,t) \\
\frac{\partial w}{\partial x}(0,t)=0,&&\frac{\partial w}{\partial x}(1,t)=\beta  u(t)\\
w(x,0)=w_0(x),&& \frac{\partial w}{\partial t}(x,0)=w_1(x,0)
\eean
Now it is convenient to have the control act at $x=1$.

Let $z(x,t)$ and ${\cal A}$ be as before (\ref{z}, \ref{mdo}).  The new boundary conditions on 
 $z(x,t)$ are
 \bean
\frac{\partial z_1}{\partial x} (0,t)=0, &\quad\ &\frac{\partial z_1}{\partial x} (1,t)=\beta u(t)
 \eean

 The open loop eigenvalues of ${\cal A}$  under Neumann boundary conditions are 
 \bea \label{laN}
  \lambda_n &=&{-\alpha+\mbox{ sign}(n) \sqrt{\alpha^2-4 n^2\pi^2}\over 2}
  \eea
  for $n=0,\pm1,\pm 2,\ldots $.
The coresponding  eigenvectors are
\bean
\phi^0(x)=\bmt 1\\0\emt,&\quad\ &\phi^n(x)=\bmt {1\over \lambda_n} \\1\emt \cos n\pi x
\eean 
for $n=\pm1,\pm 2,\ldots $.

Again we wish to stabilize this system to $z(x,t)=0$ which implies $w(x,t)=w^*(x,t)$
so we consider the Linear Quadratic Regulator problem of minimizing (\ref{crit})
with $Q(x_1,x_2)$ and $R$ as before. 
Let $P(x_1,x_2)$ be a  $2\times  2$ dimensional symmetric matrix valued function
\bean
P(x_1,x_2)&=& \bmt P_{1,1}(x_1,x_2) &P_{1,2}(x_1,x_2)\\P_{2,1}(x_1,x_2)&P_{2,2}(x_1,x_2)\emt
\eean
 which is also symmetric in $(x_1,x_2)$,
$
P(x_1,x_2)=P(x_2,x_1)
$.
Suppose there exists a control trajectory $u(t)$ such that the resulting state trajectory $z(x,t) \to 0$ as $t \to \infty$.  Then by the Fundamental 
Theorem of Calculus
\bean
0&=&\iint_{\cal S}  z'(x_1,1) P(x_1,x_2) z(x_2,0)\ dA\\
&& +\int_0^\infty \iint_{\cal S} {d\over dt}\left( z'(x_1,t) P(x_1,x_2) z(x_2,t)\right)\ dA \ dt
\eean

Now we assume that $P(x_1,x_2)$  satisfies Neumann boundary conditions in both its arguments
\bean
\frac{\partial P}{\partial x_1}(0,x_2)&=& \frac{\partial P}{\partial x_1}(1,x_2)\ =\ 0\\
\frac{\partial P}{\partial x_2}(x_1,0)&=& \frac{\partial P}{\partial x_2}(x_1,1)\ =\ 0
\eean
and integrate by parts as before to obtain
\bean
0&=&\iint_{\cal S}  z'(x_1,1) P(x_1,x_2) z(x_2,0)\ dA\\
&& +\int_0^\infty \iint_{\cal S} z_2(x_1,t)P_{1,1}(x_1,x_2)z_1(x_2,t)+ z_1(x_1,t)P_{1,1}(x_1,x_2)z_2(x_2,t)\\
&&+z_2(x_1,t)P_{1,2}(x_1,x_2)z_2(x_2,t)+ z_2(x_1,t)P_{2,1}(x_1,x_2)z_2(x_2)
\\&& -\alpha z_2(x_1,t)P_{2,1}(x_1,x_2)z_1(x_2)-\alpha z_1(x_1,t)P_{1,2}(x_1,x_2)z_2(x_2)
\\ && -2\alpha z_2(x_1,t)P_{2,2}(x_1,x_2)z_2(x_2)\\
&&+  z_1(x_1,t)\frac{\partial^2 P_{2,1}}{\partial x_1^2}(x_1,x_2)z_1(x_2,t)+z_1(x_1,t)\frac{\partial^2 P_{1,2}}{\partial x_2^2}(x_1,x_2)z_1(x_2,t)\\
&&+  z_1(x_1,t)\frac{\partial^2 P_{2,2}}{\partial x_1^2}(x_1,x_2)z_2(x_2,t)+z_2(x_1,t)\frac{\partial^2 P_{2,2}}{\partial x_2^2}(x_1,x_2)z_1(x_2,t)
\ dA \ dt \\
&&+\int_0^\infty  \int_0^1 \beta u(t) P_{2,1} (1,x_2) z_1(x_2,t)\ dx_2\ dt
\\
&&+\int_0^\infty  \int_0^1 z_1(x_1,t)P_{1,2}(x_1,1)\beta u(t) \ dx_1\ dt
\\
&&+\int_0^\infty  \int_0^1 \beta u(t)P_{2,2}(1,x_2) z_2(x_2,t)\ dx_2\ dt
\\
&&+\int_0^\infty  \int_0^1 z_2(x_1,t)P_{2,2}(x_1,1)\beta u(t) \ dx_1\ dt
\eean

We add this to the criterion (\ref{crit}) to get the equivalent criterion to be minimized
\bean
&&\int_0^\infty \iint_{\cal S}  z'(x_1,t) Q(x_1,x_2) z(x_2,t)\ dA+R( u(t))^2\ dt\\
&&+\iint_{\cal S}  z'(x_1,1) P(x_1,x_2) z(x_2,0)\ dA\\
&& +\int_0^\infty \iint_{\cal S} z_2(x_1,t)P_{1,1}(x_1,x_2)z_1(x_2,t)+ z_1(x_1,t)P_{1,1}(x_1,x_2)z_2(x_2,t)\\
&&+z_2(x_1,t)P_{1,2}(x_1,x_2)z_2(x_2,t)+ z_2(x_1,t)P_{2,1}(x_1,x_2)z_2(x_2)
\\&& -\alpha z_2(x_1,t)P_{2,1}(x_1,x_2)z_1(x_2)-\alpha z_1(x_1,t)P_{1,2}(x_1,x_2)z_2(x_2)
\\ && -2\alpha z_2(x_1,t)P_{2,2}(x_1,x_2)z_2(x_2)\\
&&+  z_1(x_1,t)\frac{\partial^2 P_{2,1}}{\partial x_1^2}(x_1,x_2)z_1(x_2,t)+z_1(x_1,t)\frac{\partial^2 P_{1,2}}{\partial x_2^2}(x_1,x_2)z_1(x_2,t)\\
&&+  z_1(x_1,t)\frac{\partial^2 P_{2,2}}{\partial x_1^2}(x_1,x_2)z_2(x_2,t)+z_2(x_1,t)\frac{\partial^2 P_{2,2}}{\partial x_2^2}(x_1,x_2)z_1(x_2,t)
\ dA \ dt \\
&&+\int_0^\infty  \int_0^1 \beta u(t) P_{2,1} (1,x_2) z_1(x_2,t)\ dx_2\ dt
\\
&&+\int_0^\infty  \int_0^1 z_1(x_1,t)P_{1,2}(x_1,1)\beta u(t) \ dx_1\ dt
\\
&&+\int_0^\infty  \int_0^1 \beta u(t)P_{2,2}(1,x_2) z_2(x_2,t)\ dx_2\ dt
\\
&&+\int_0^\infty  \int_0^1 z_2(x_1,t)P_{2,2}(x_1,1)\beta u(t) \ dx_1\ dt
\eean

We complete the square as before 
to obtain 
\bea \label{K1N}
K(x_2)&=&\bmt K_1(x_2)&K_2(x_2)\emt \ =\  -R^{-1} \beta  \bmt P_{2,1} (1,x_2) & P_{2,2} (1,x_2) \emt 
\eea
and
\bea
\nonumber&&\frac{\partial^2 P_{2,1}}{\partial x_1^2}(x_1,x_2) + \frac{\partial^2 P_{1,2}}{\partial x_2^2}(x_1,x_2)+Q_{1,1}(x_1,x_2)\\
&&=\gamma^2P_{1,2}(x_1,1)P_{2,1}(1,x_2)
\label{RPDE11aN}\\
\nonumber&&P_{11}(x_1,x_2) -\alpha P_{1,2}(x_1,x_2) +\frac{\partial^2 P_{2,2}}{\partial x_1^2}(x_1,x_2)+Q_{1,2}(x_1,x_2)   
\\
&&=\gamma^2P_{1,2}(x_1,1)P_{2,2}(1,x_2)
\label{RPDE12aN}\\
\nonumber&&P_{11}(x_1,x_2)-\alpha P_{2,1}(x_1,x_2)+\frac{\partial^2 P_{2,2}}{\partial x_2^2}(x_1,x_2)+Q_{21}(x_1,x_2) \\
&&=\gamma^2P_{2,2}(x_1,1)P_{2,1}(1,x_2)
 \label{RPDE21aN}\\
\nonumber&&P_{1,2}(x_1,x_2)+P_{2,1}(x_1,x_2)-2\alpha P_{2,2}(x_1,x_2)+Q_{22}(x_1,x_2)\\
&&=\gamma^2P_{2,2}(x_1,1)P_{2,1}(1,x_2)
\label{RPDE22aN}
\eea
where $\gamma^2=R^{-1}\beta^2$.
We call these equations (\ref{RPDE11aN},\ref{RPDE12aN},\ref{RPDE21aN},\ref{RPDE22aN}) the Riccati PDE 
for the LQR Neumann boundary  control of the wave equation.

Since $P(x_1,x_2)$ satisfies Neumann boundary conditions   at $x=0$ and $x=1$ it can be expressed in terms of the eigenfunctions of $\frac{\partial^2}{\partial x^2}$
with respect to Neumann  boundary conditions.  As before it is convenient to decouple the spatial frequencies.
\bea  \label{PsumN}
P(x_1,x_2)&=&  \sum_{n=0}^\infty \bmt P^{n,n}_{1,1}&P^{n,n}_{1,2}\\P^{n,n}_{2,1}&P^{n,n}_{2,2}\emt \cos n\pi x_1  \cos n\pi x_2
\eea
and  $P^{n,n}_{1,2}=P^{n,n}_{2,1}$.
Then (\ref{K1N}) implies 
\bea \label{KsumN}
K(x_2)&=&-R^{-1}\beta    \sum_{n=0}^\infty \bmt P_{2,1}^{n,n} &P_{2,2}^{n,n} \emt  \cos n\pi x 
\eea
 
We assume that $Q(x_1,x_2)$ has a similar expansion
\bea \label{QsumN}
Q(x_1,x_2)&=&  \sum_{n=0}^\infty  Q^{n,n} \cos n\pi x_1  \cos n\pi x_2
\eea
This leads to an infinite algebraic Riccati equation for $P^{n,n}$.

The Riccati PDE (\ref{RPDE11aN},\ref{RPDE12aN},\ref{RPDE21aN},\ref{RPDE22aN}) implies that 
\bea
0&=&   -2n^2\pi^2 P^{n,n}_{1,2} +Q^{n,n}_{1,1}-\gamma^2\left(P^{n,n}_{1,2}\right)^2  \label{11N}\\
0&=& P^{n,n}_{1,1}-\alpha P^{n,n}_{1,2} -n^2\pi^2 P^{n,n}_{2,2}+Q^{n,n}_{1,2} -\gamma^2P^{n,n}_{1,2}P^{n,n}_{2,2} \label{12N}\\
0&=& P^{n,n}_{1,1} -\alpha P^{n,n}_{2,1}-n^2\pi^2 P^{n,n}_{2,2}+Q^{n,n}_{2,1} -\gamma^2P^{n,n}_{2,2}P^{n,n}_{2,1}\label{21N}\\
0&=& 2P^{n,n}_{1,2}-2\alpha P^{n,n}_{2,2}+Q^{n,n}_{2,2}- \gamma^2\left(P^{n,n}_{2,2}\right)^2 \label{22N}
\eea
where $\gamma^2= R^{-1}\beta^2$.

For each $n$ these are the Riccati equations of  the two dimensional  LQR with matrices 
\eq  \label{lqr2N}
F^{n,n}=\bmt 0&1\\ -n^2 \pi^2 &-\alpha\emt, && G^{n,n}=\bmt 0\\ \beta\emt\\
Q^{n,n}=\bmt Q^{n,n}_{1,1}& Q^{n,n}_{1,2}\\Q^{n,n}_{2,1}&Q^{n,n}_{2,2}\emt,&& R^{n,n}=\bmt R\emt
\ee
Notice $G^{n,n}$ is different from before.

We solve these equations and obtain 
\bea 
P^{n,n}_{1,1}&=& \alpha P^{n,n}_{1,2}+\left(n^2\pi^2 +\gamma^2 P^{n,n}_{1,2}  \right) P^{n,n}_{2,2}-Q^{n,n}_{1,2}  \label{p11N}\\
\label{p12N}
P^{n,n}_{1,2}&=& {-n^2\pi^2\pm \sqrt{n^4\pi^4+Q^{n,n}_{1,1}}\over \gamma^2}\\
P^{n,n}_{2,2}&=&{-\alpha \pm \sqrt{\alpha^2+\gamma^2\left(Q^{n,n}_{2,2}+2P^{n,n}_{1,2}\right)}\over \gamma^2} \label{p22N}
\eea
Since we want $P^{n,n}_{2,2}$ to be nonegative we need to take the positive sign in (\ref{p22N}).
If the two dimensional LQR (\ref{lqr2N}) satisfies the standard conditions then the associated
 Riccati equation has a unique nonnegative definite solution.  This implies that if we take the
  negative sign in (\ref{p12N}) the resulting $P^{n,n}$ is not nonnegative definite.

The $2\times 2$ closed loop system is 
\bean
F^{n,n}+G^{n,n}K^{n,n}&=& \bmt 0&1\\ -n^2\pi^2 -\gamma^2 P_{2,1}^{n,n}&-\alpha -\gamma^2 P_{2,2}^{n,n}\emt
\eean
and the closed loop eigenvalues are
\bea\label{clevaN}
\mu_n&=& 
-{\alpha +\gamma^2P^{n,n}_{2,2}\over 2}+\left(\mbox{sign n}\right){\sqrt{(\alpha + \gamma^2P^{n,n}_{2,2})^2-4(n^2\pi^2+\gamma^2 P^{n,n}_{2,1})}\over 2}
\eea
for $n=0,\pm 1,\pm2,\ldots$.
For $n=0,1,2,\ldots$ the corresponding eigenvectors of the $n^{th}$ $2\times 2$ closed loop system are 
 \bea \label{eigv2N}
 \bmt {1\over \mu_n}\\1 \emt,& \quad\ & \bmt {1\over \mu_{-n}}\\1 \emt
 \eea
 The corresponding eigenvectors of the infinite closed loop system are 
 \bea \label{eigv2inN}
v_n(x)= \bmt {1\over \mu_n}\\ 1\emt \cos n\pi x ,& \quad\ & \bmt {1\over \mu_{-n}}\\ 1\emt \cos n\pi x
 \eea
The trajectories of the infinite closed loop system are 
\bea
z(x,t)&=& \sum_{n=-\infty}^\infty \zeta_n(t) \bmt {1\over \mu_n}\\1 \emt \cos n\pi x
\eea
where
\bean
\zeta_n(t) &=& e^{\mu_n t}  \zeta^0
\eean

The  quadratic kernel of the optimal cost is given by (\ref{Psum}) which reduces to 
\bea \label{Psum_nN}
P(x_1,x_2)&=& \sum_{n=1}^\infty \bmt P_{1,1}^{n,n}&P_{1,2}^{n,n}\\P_{2,1}^{n,n}&P_{2,2}^{n,n}\emt \cos n\pi x_1 \cos n\pi x_2
\eea
Again we can control each spatial frequency indepedently.   

But again we face the question  is (\ref{Psum_nN}) convergent.
And if the system is undamped, $\alpha=0$,
  how many spatial modes can we stabilize and how fast?  To partially answer these questions we look at our  simple example again,  we assume 
\bean 
\beta=1\\
R=1\\
Q^{0,0}_{1,1}=Q^{0,0}_{2,2}=q\\
Q^{n,n}_{1,1}=Q^{n,n}_{2,2}={q\over n^r }  \mbox{ for } n>0\\
Q^{n,n}_{1,2}=Q^{n,n}_{2,1}=0 \mbox{ for } n\ge 0
\eean
For  (\ref{QsumN}) to converge we must take $r>1$. 


For this  system
\bean
\gamma^2&=&1\\
P^{n,n}_{1,2}&=&-n^2\pi^2+ \sqrt{n^4\pi^4+{q\over n^{r}}}\\
P^{n,n}_{2,2}&=&\sqrt{q/n^r+2P^{n,n}_{1,2}}
\\
P^{n,n}_{1,1}&=&\left( n^2\pi^2+P^{n,n}_{1,2}\right) P^{n,n}_{2,2}
\eean

Again  
by the Mean Value Theorem there exists an $s$ between $ n^2\pi^2$ and  $n^2\pi^2+{q\over n^r}$
such that
\bean
P^{n,n}_{1,2}&=& \sqrt{n^4\pi^4+{q\over n^{r}}} -n^2\pi^2\ =\  {q\over 2s^{1/2}n^r}
\eean 
The function ${1\over 2s^{1/2}n^r}$ is montonically decreasing between $ n^2\pi^2$ and  $n^2\pi^2+{q\over n^r}$
and takes on its maximum value at 	 $ s=n^2\pi^2$ so
\bean
0&<&P^{n,n}_{1,2}\le   {q\over 2\pi n^{r+1}}
\eean
Then
\bean
P^{n,n}_{2,2}&=&\sqrt{q/n^r+2P^{n,n}_{1,2}} \ =\ O(n^{-r/2})
\\
P^{n,n}_{1,1}&= &\left( n^2\pi^2+P^{n,n}_{1,2}\right) P^{n,n}_{2,2}\ = \ O(n^{2-r/2})
\eean
This time for  $P_{1,1}(x_1,x_2)$ in the optimal cost  (\ref{PsumN}) to converge  we must take $r>6$ so $2-r/2<-1$.  
Then $P^{n,n}_{2,2}=O(n^{-r/2})$ so the higher spatial modes are less damped.   Moreover the first term inside the the square root in (\ref{clevaN}) is going to zero while the negative term is   growing so at least the higher eigenvalues are complex.

But for the optimal feedback (\ref{KsumN}) to converge  we only need to take $r>2$. If we take $r$ a little larger than $2$
then the term outside the square root in the closed loop eigenvalues (\ref{clevaN}) is decaying  so the higher spatial modes are less damped.   Hence we conclude for this simple example Dirichlet boundary control of the wave equation is preferable to Neumann boundary control.  We strongly suspect that this is true in general.

\section{Conclusion}
We have used the simple and constructive technique of completing the square to solve the LQR problems for the linear wave equation under both Dirichlet and Neumann boundary control.  The results are explicit formulas for the quadratic optimal cost and the linear optimal feedback.  These formulas decouple the spatial frequencies so we can damp out all or just some frequencies.   We can also use the linear optimal feedback to stabilize the wave to any desired  trajectory 
of the open loop system.   We would like to thank Miroslav Krstic for his helpful comments.

\end{document}